\documentclass[11pt,a4paper]{article}
\usepackage{amsmath,amsthm,amsfonts,amssymb,mathrsfs}

\usepackage{graphicx}
\usepackage{subfigure}

\usepackage[latin1]{inputenc}
\usepackage[swedish,english]{babel}

\begin{document}

\title{A higher-order singularity subtraction technique for the
  discretization of singular integral operators on curved
  surfaces\thanks{Supported by the Swedish Research Council under
    contract 621-2011-5516.}}
\author{Johan Helsing\\
  Centre for Mathematical Sciences\\
  Lund University, P.O. Box 118, SE-221 00 Lund, Sweden}

\maketitle


\section{Introduction}

A major challenge facing the community working on integral equation
based solvers for boundary value problems is the construction of
efficient discretizations of integral operators with singular kernels
on curved surfaces. Classic approaches such as singularity
subtraction, special purpose quadrature, singularity cancellation,
kernel regularization, and various adaptive strategies may work well
in many situations but have not yet fully, in three dimensions,
succeeded in unleashing the computational power of integral equation
methods needed for excellence in real-world physics
applications~\cite[Section 1]{Klock12}. Recently, two new promising
methods have been launched: the quadrature by extension (QBX) method
which exploits that fields induced by integral operators are often
smooth close to the boundaries where their sources are
located~\cite{Klock12} and a method relying on a combination of
adaptivity, local invertible affine mappings with certain
orthogonality properties, and the use of precomputed tables of
quadrature rules~\cite{Brem12}. It seems to be an open question what
method, or combination of techniques, is best.

This note is about promoting a classic technique for the
discretization of singular integral operators on curved surfaces,
namely singularity subtraction. The idea is to use analytical
evaluation to a maximum degree and split singular (and nearly
singular) operators into two parts each -- one ill-behaved part whose
action can be evaluated using high-order analytic product integration,
and another more regular part for which purely numerical integration
is used, compare~\cite{Farin01}. Based on this idea we present and
implement a simple Nystr{\"o}m scheme for Laplace's equation on tori.
Surprisingly accurate results are produced.

\section{Problem formulation and methods}

We consider the interior Dirichlet Laplace problem
\begin{align}
\Delta U(r)&=0\,, \qquad\quad\; r\in V\,,\label{eq:PDE1}\\
       U(r)&=g(r)\,, \qquad r\in \Gamma\,,\label{eq:PDE2}
\end{align}
where $g(r)$ is a smooth function on the boundary $\Gamma$ of a smooth
domain $V$ in $\mathbb{R}^3$. For the solution
of~(\ref{eq:PDE1},\ref{eq:PDE2}) we use the double-layer
representation
\begin{equation}
U(r)=\frac{1}{4\pi}\int\int_\Gamma\frac{n_{r'}\cdot(r'-r)}{|r'-r|^3}
\mu(r')\,{\rm d}\sigma_{r'}\,,
\label{eq:rep}
\end{equation}
where $n_r$ is the exterior unit normal of $\Gamma$ at position $r$,
${\rm d}\sigma$ is an element of surface area, and $\mu$ is an unknown
layer density. An integral equation formulation
for~(\ref{eq:PDE1},\ref{eq:PDE2}) reads
\begin{equation}
\mu(r)+\frac{1}{2\pi}\int\int_\Gamma\frac{n_{r'}\cdot(r'-r)}{|r'-r|^3}
\mu(r')\,{\rm d}\sigma_{r'}=2g(r)\,.
\label{eq:inteq}
\end{equation}

\subsection{Parameterization}

The domain $V$ is taken to be a torus whose surface $\Gamma$ is
parameterized over the square $\left\{s=(s_1,s_2)\in\mathbb{R}^2:
  -\pi\le s_1,s_2\le\pi\right\}$ as
\begin{equation}
r(s)=\left[\varrho(s)\cos(s_2), \varrho(s)\sin(s_2), 
\delta_2\sin(s_1) \right]\,,
\label{eq:param}
\end{equation}
where
\begin{equation}
\varrho(s)=2+\delta_1\cos(2s_2)+\delta_2\cos(s_1)
\end{equation}
and $\delta_1$ and $\delta_2$ are shape parameters. The choice
$\delta_1=0$ corresponds to the standard tori used in~\cite[Section
3.5]{Brem12}.

\begin{figure}
\centering
\includegraphics[height=32mm]{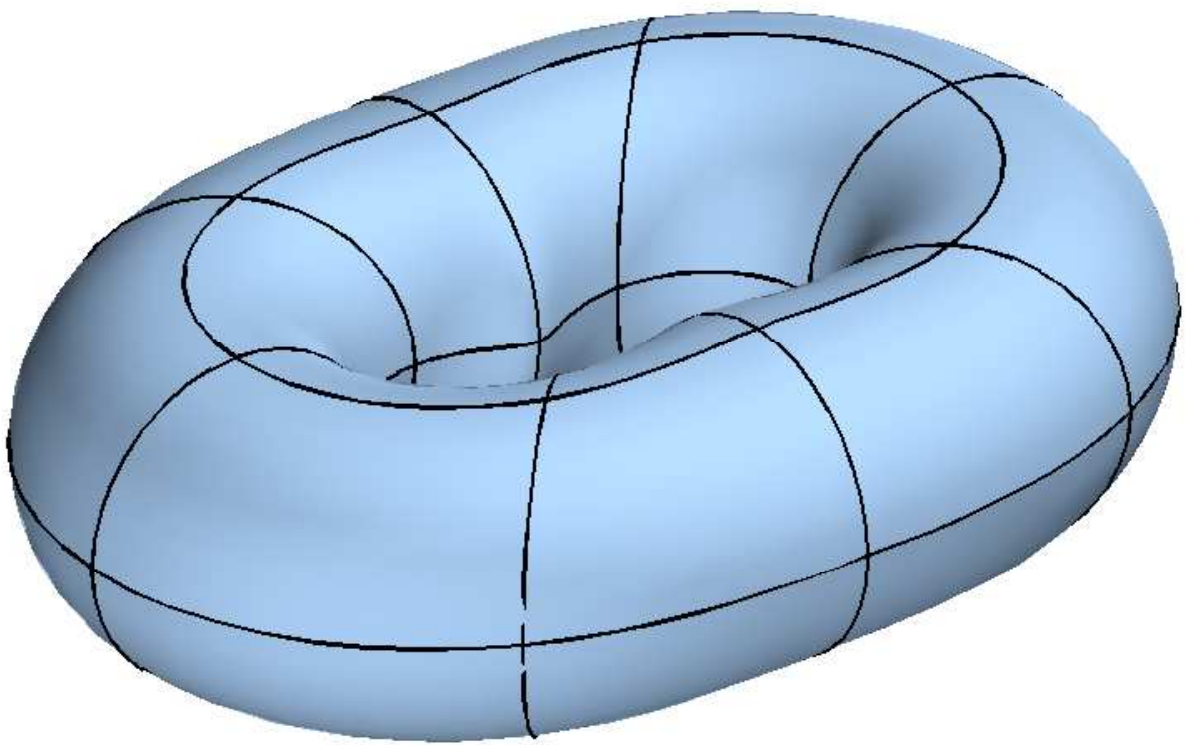}
\hspace{8mm}
\includegraphics[height=32mm]{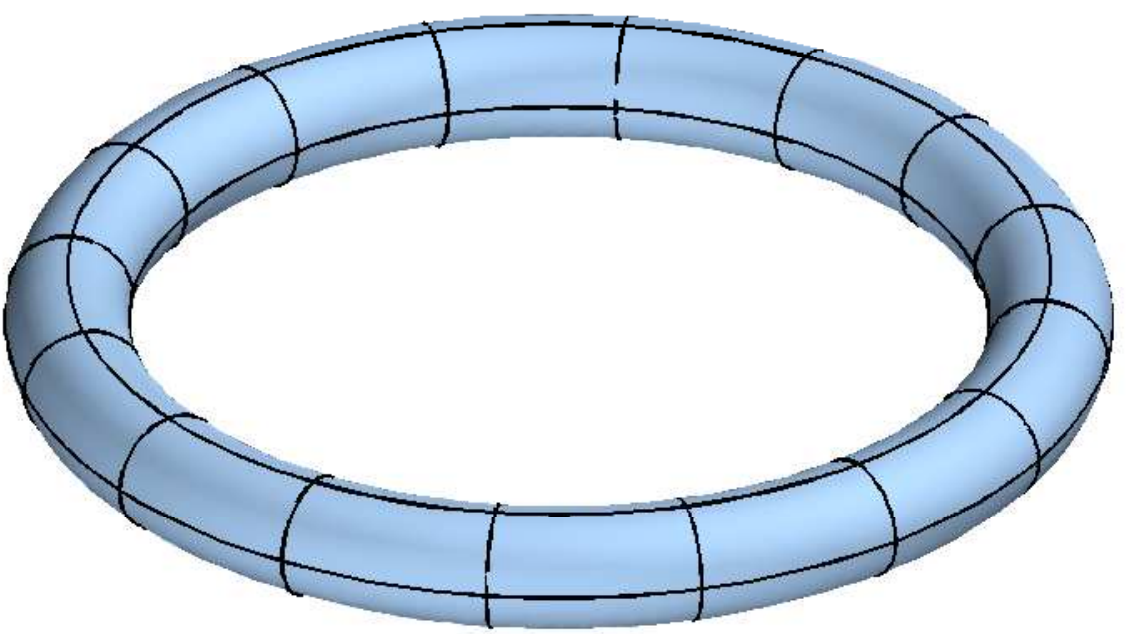}
\caption{\sf Left: a torus with $\delta_1=0.5$,
  $\delta_2=1$, $p_1=4$, $p_2=8$ and 32 patches $\Gamma_{ij}$. Right:
  $\delta_1=0$, $\delta_2=0.25$, $p_1=4$, $p_2=16$ and 64 patches
  $\Gamma_{ij}$.}
\label{fig:tori}
\end{figure}

We shall use Nystr{\"o}m discretization for~(\ref{eq:inteq}) based on
composite tensor product Gauss--Legendre quadrature. For this, we
introduce a sequence of mappings $\rho_{ij}$ with $i=1,\ldots,p_1$ and
$j=1,\ldots,p_2$
\begin{equation}
\rho_{ij}(t)=r(\pi(t_1+2i-p_1-1)/p_1,\pi(t_2+2j-p_2-1)/p_2)\,.
\label{eq:param2}
\end{equation}
The mapping $\rho_{ij}(t)$ covers a patch $\Gamma_{ij}$ of $\Gamma$
when mapped from the square $\left\{t=(t_1,t_2)\in\mathbb{R}^2: -1\le
  t_1,t_2\le 1\right\}$. The disjoint union of the $\Gamma_{ij}$ is
$\Gamma$. See Figure~\ref{fig:tori} for two examples.

Using~(\ref{eq:param2}) and introducing
\begin{align}
u_{ij}(r,t')&=\rho_{ij}(t')-r\,,\\
J_{ij}(r,t')&=\left(\frac{\partial\rho_{ij}(t')}{\partial t'_2}\times
  \frac{\partial \rho_{ij}(t')}{\partial t'_1}\right)\cdot u_{ij}(r,t')\,,
\end{align}
we can rewrite the integral operator in~(\ref{eq:inteq}) as the sum
\begin{equation}
\sum_{i=1,j=1}^{p_1,p_2}D_{ij}\mu(r)=
\sum_{i=1,j=1}^{p_1,p_2}\frac{1}{2\pi}\int_{-1}^1\int_{-1}^1
\frac{J_{ij}(r,t')}{|u_{ij}(r,t')|^3}\mu(r')\,{\rm d}t'_1\,{\rm d}t'_2\,.
\label{eq:Dsum}
\end{equation}

\subsection{Singularity subtraction}

Nystr{\"o}m discretization works well for a particular $D_{ij}$
in~(\ref{eq:Dsum}) if $r$ is far away from the patch $\Gamma_{ij}$.
If $r$ is close to, or on, $\Gamma_{ij}$ then the kernel is nearly
singular, or singular, and something better is needed. Let
\begin{align}
v_{ij}(r,t')&=
\sum_{k=1}^2(t'_k-t_k(r))\frac{\partial\rho_{ij}}{\partial
  t_k}(t(r))\,,\label{eq:vij}\\
\Delta_{ij}(r,t')&=|u_{ij}(r,t')|^2-|v_{ij}(r,t')|^2\,,
\end{align}
where $t(r)=\rho_{ij}^{-1}(r)$. For $t'$ close to $t(r)$ the operator
$D_{ij}$ can be expanded
\begin{equation}
D_{ij}=\sum_{k=0}^{\infty}D_{ijk}\,,
\label{eq:expand}
\end{equation}
\begin{equation}
D_{ijk}\mu(r)=\binom{-3/2}{k}\frac{1}{2\pi}\int_{-1}^1\int_{-1}^1
\frac{J_{ij}(r,t')\Delta_{ij}^k(r,t')}{|v_{ij}(r,t')|^{3+2k}}
\mu(r')\,{\rm d}t'_1\,{\rm d}t'_2\,.
\label{eq:Dijk}
\end{equation}
See~\cite[Section 1]{Johns89} for a discussion of a similar expansion.

Our proposed singularity subtraction technique for $t'$ close to
$t(r)$ makes use of the split
\begin{equation}
D_{ij}=D_{ij}^K+D_{ij}^{\circ}\,.
\label{eq:split}
\end{equation}
Here
\begin{equation}
D_{ij}^K=\sum_{k=0}^KD_{ijk}\,, \qquad 
D_{ij}^{\circ}=D_{ij}-D_{ij}^K\,,
\label{eq:split2}
\end{equation}
with $K$ a small integer and with $D_{ij}$ as in~(\ref{eq:Dsum}). The
action of $D_{ij}^K$ is to be evaluated using high-order analytic
product integration and $D_{ij}^{\circ}$ is supposed to be
sufficiently smooth as to allow for accurate Nystr{\"o}m
discretization.

\section{Recursive evaluation of integrals}
\label{sec:recur}

Computing $D_{ijk}\mu(r)$ in~(\ref{eq:Dijk}) requires the evaluation
of an expression of the form
\begin{equation}
\sum_{m,n=0}\int_{-1}^1\int_{-1}^1\frac{\alpha_{mn}{t'}_1^m {t'}_2^n
\,{\rm d}t'_1\,{\rm d}t'_2}
{\left(a^2(t'_1-t_1)^2+2abc(t'_1-t_1)(t'_2-t_2)+b^2(t'_2-t_2)^2\right)^{k+3/2}}
\,,
\label{eq:form1}
\end{equation}
where $a$, $b$, and $c$ are constants known from~(\ref{eq:vij}) and
$\alpha_{mn}$ are coefficients of a polynomial approximating the
smooth function $J_{ij}(r,t')\Delta_{ij}^k(r,t')\mu(r')$. The variable
substitution $x=at'_1$, $x_0=at_1$, $y=bt'_2$, $y_0=bt_2$ makes the
terms in~(\ref{eq:form1}) appear as
\begin{equation}
\frac{1}{a^{m+1}b^{n+1}}
\int_{-a}^a\int_{-b}^b\frac{\alpha_{mn}x^m y^n\,{\rm d}x\,{\rm d}y}
{\left((x-x_0)^2+2c(x-x_0)(y-y_0)+(y-y_0)^2\right)^{k+3/2}}\,.
\label{eq:form2}
\end{equation}

We now present a scheme for the evaluation of integrals of the
form~(\ref{eq:form2}). Let
\begin{equation}
d_c(x,y)=x^2+2cxy+y^2
\end{equation}
and define in the Hadamard finite part sense the indefinite integrals
\begin{align}
C_{mnk}(x,x_0,y,y_0,c)&=\int\int\frac{x^m y^n\,{\rm d}x\,{\rm d}y}
{|d_c(x-x_0,y-y_0)|^{k+1/2}}\,,\\
F_{mk}(x,x_0,y,y_0,c)&=\int\frac{x^m\,{\rm d}x}
{|d_c(x-x_0,y-y_0)|^{k+1/2}}\,,\\
G_{nk}(x,x_0,y,y_0,c)&=\int\frac{y^n\,{\rm d}y}
{|d_c(x-x_0,y-y_0)|^{k+1/2}}\,.
\end{align}
Using partial integration, and for given values $(x,x_0,y,y_0,c)$, one
can show that when $m+n+1\ne 2k$ holds
\begin{equation}
C_{mnk}=\frac
{mx_0C_{(m-1)nk}+ny_0C_{m(n-1)k}+(x-x_0)x^mG_{nk}+(y-y_0)y^nF_{mk}}
{m+n+1-2k}\,.
\label{eq:recur1}
\end{equation}
When $m+n+1=2k$ holds
\begin{align}
C_{mnk}&=x_0C_{(m-1)nk}+\beta_k\Big((m-1)C_{(m-2)n(k-1)}-cnC_{(m-1)(n-1)(k-1)}
\nonumber\\
&\quad-x^{m-1}G_{n(k-1)}+cy^nF_{(m-1)(k-1)}\Big)\,,\label{eq:recur2}\\
C_{mnk}&=y_0C_{m(n-1)k}+\beta_k\Big((n-1)C_{m(n-2)(k-1)}-cmC_{(m-1)(n-1)(k-1)}
\nonumber\\
&\quad-y^{n-1}F_{m(k-1)}+cx^mG_{(n-1)(k-1)}\Big)\,,\label{eq:recur3}
\end{align}
where $\beta_k=1/((1-c^2)(2k-1))$. For the evaluation of $F_{mk}$ we
use
\begin{align}
F_{00}&=\log\left(|d_c(x-x_0,y-y_0)|^{1/2}+(x-x_0)+c(y-y_0)\right)\,,
\label{eq:recur4}\\
F_{0k}&=\frac{\beta_k}{(y-y_0)^2}
\left(\frac{(x-x_0)+c(y-y_0)}{|d_c(x-x_0,y-y_0)|^{k-1/2}}
+2(k-1)F_{0(k-1)}\right)\,,\quad k\ge 1\,,\label{eq:recur5}\\
F_{m0}&=\frac{1}{m}\Big( 
x^{m-1}{|d_c(x-x_0,y-y_0)|^{1/2}}-(m-1)d_c(-x_0,y-y_0)F_{(m-2)0}\nonumber\\
&\quad+(2m-1)(x_0-c(y-y_0))F_{(m-1)0}\Big)\,,\quad m\ge 1\,,\label{eq:recur6}\\
F_{mk}&=\frac{1}{2k-1}\left((m-1)F_{(m-2)(k-1)}
-x^{m-1}{|d_c(x-x_0,y-y_0)|^{k-1/2}}\right)\nonumber\\
&\quad+(x_0-c(y-y_0))F_{(m-1)k}\,,\quad m,k\ge 1\,.
\label{eq:recur7}
\end{align}
Expression for $G_{nk}$ are obtained by interchanging
$x\rightleftarrows y$ and $m\rightleftarrows n$ in the expressions for
$F_{mk}$. The recursions for $C_{mnk}$, $F_{mk}$ and $G_{nk}$ allow
the integrals in~(\ref{eq:form1}) to be evaluated at a modest
computational cost. Note that, for this, only $C_{mnk}$ with $m,n\ge
0$ and $k\ge 1$ are needed.

\section{Details on the discretization}
\label{sec:disc}

We now give precise details on our Nystr{\"o}m discretization
of~(\ref{eq:inteq}). Aiming at 10th order convergence we take 10-point
composite tensor product Gauss--Legendre quadrature (GL10) as our
underlying quadrature scheme. On each $\Gamma_{ij}$ there will then be
a grid of $100$ discretization points where the discretized density
$\boldsymbol{\mu}$ is sought. The discretized system~(\ref{eq:inteq})
has $100p_1p_2$ unknowns. We shall also use a temporary, finer, grid
with $256$ discretization points on each $\Gamma_{ij}$ placed
according to 16-point composite tensor product Gauss--Legendre
quadrature (GL16).

If, for a particular $D_{ij}$ and $r$ in~(\ref{eq:Dsum}), the local
parameter $t=\rho_{ij}^{-1}(r)$ is such that $3.5<|t|$, then the point
$r$ is considered far away from $\Gamma_{ij}$ and we discretize
$D_{ij}\mu(r)$ using the underlying GL10 scheme.

If $2<|t|\le 3.5$, then $r$ is somewhat close to $\Gamma_{ij}$ and we
use an extended scheme: first $\boldsymbol{\mu}$ is interpolated to
the finer grid on $\Gamma_{ij}$ and then $D_{ij}\mu(r)$ is discretized
using GL16. High-degree polynomial interpolation of smooth functions
known at Legendre nodes can be very accurate, despite involving
ill-conditioned Vandermonde systems~\cite[Appendix A]{Hels08}.

If $|t|<2$, then $r$ is close to $\Gamma_{ij}$, the operator $D_{ij}$
is nearly singular or weakly singular, and we use the
split~(\ref{eq:split}). The discretization is carried out on the GL16
grid on $\Gamma_{ij}$, which means that $\boldsymbol{\mu}$ has to be
interpolated to 256 points as in the previous paragraph. The operator
$D_{ij}^{\circ}$ is discretized using GL16. The operator $D_{ij}^K$ is
discretized using the method of Section~\ref{sec:recur}. We let
$m,n=0,\ldots,15$ in~(\ref{eq:form1}). The 256 coefficients
$\alpha_{mn}$ are obtained by multiplying the pointwise values of
$\boldsymbol{\mu}$ at the 256 fine grid points on $\Gamma_{ij}$ with
pointwise values of $J_{ij}(r,t')\Delta_{ij}^k(r,t')$ and then, in
principle, solving a Vandermonde system of size $256\times 256$. In
practice one can obtain the $\alpha_{mn}$ by solving two $16\times 16$
systems with multiple right hand sides. As for the optimal number $K$
in~(\ref{eq:split}), it turns out to be related to the polynomial
degree of the underlying discretization and to the overall mesh
refinement determined by $p_1$ and $p_2$ of~(\ref{eq:param2}). For
small vales $p_1$ and $p_2$ and high degree quadrature, $K$ should be
rather low.  We choose $K=1$, that is, we use two terms in the sum
of~(\ref{eq:split2}).

\section{Numerical examples}

Numerical experiments are performed on tori given by~(\ref{eq:param})
using a program solely implemented in {\sc Matlab} and executed on a
workstation equipped with an IntelXeon E5430 CPU at 2.66 GHz and 32 GB
of memory. Three different $\delta=(\delta_1,\delta_2)$ are chosen:
$\delta=(0,1)$, $\delta=(0.5,1)$, and $\delta=(0,0.25)$. See
Figure~\ref{fig:tori} for illustrations. The boundary condition $g(r)$
in~(\ref{eq:PDE2}) is taken as $g(r)=1/|r-r_1|-1/|r-r_2|$, with
$r_1=(4,0,0)$ and $r_2=(0,4,0)$ for $\delta=(0,1)$, with
$r_1=(4.5,0,0)$ and $r_2=(0,3.5,0)$ for $\delta=(0.5,1)$, and with
$r_1=(3.25,0,0)$ and $r_2=(0,3.25,0)$ for $\delta=(0,0.25)$. The
discretized system~(\ref{eq:inteq}) is solved iteratively using GMRES.

\begin{figure}
\centering
\includegraphics[height=70mm]{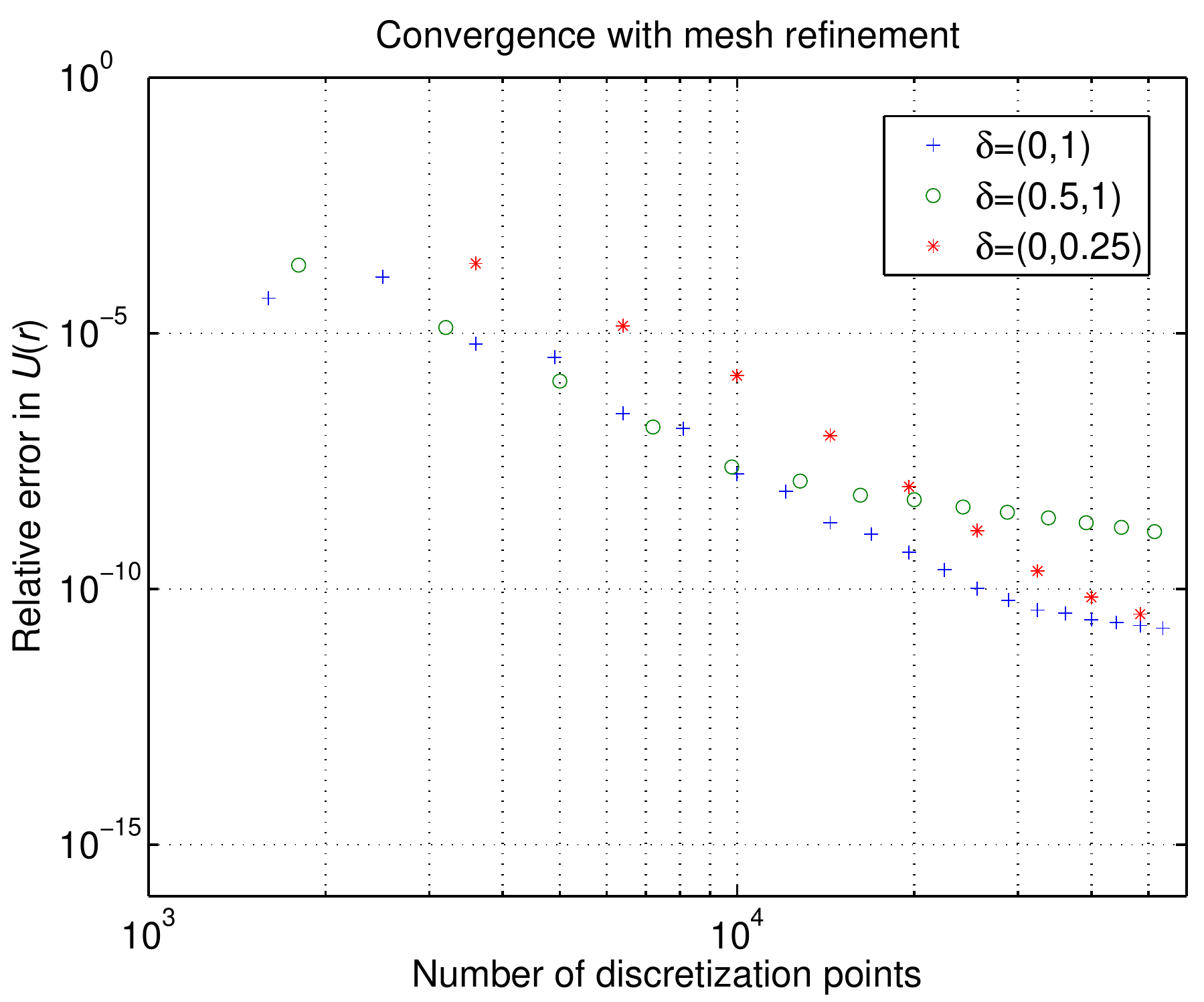}
\caption{\sf Relative $L^2$ error in $U(r)$ at the tube center when 
  solving an interior Dirichlet Laplace problem on tori given
  by~(\ref{eq:param}).}
\label{fig:conv}
\end{figure}

Figure~\ref{fig:conv} shows convergence of $U(r)$, evaluated via a
discretization of~(\ref{eq:rep}), at points along the center of the
torus tubes. The mesh is refined by increasing the parameters $p_1$
and $p_2$ of~(\ref{eq:param2}), keeping $p_2=p_1$ for $\delta=(0,1)$,
$p_2=2p_1$ for $\delta=(0.5,1)$, and $p_2=4p_1$ for $\delta=(0,0.25)$.
A relative residual less than $\epsilon_{\rm mach}$ is obtained in
between 15 and 20 iterations for reasonably resolved systems. The
recursion of Section~\ref{sec:recur} is rather fast. For example, with
10,000 discretization points and $\delta=(0,0.25)$ only 28 seconds are
spent doing singularity subtraction.

One can see in Figure~\ref{fig:conv} that the initial convergence of
$U(r)$ is approximately 10th order, as expected. As the number of
discretization points grows, however, the error stemming from the
discretization of $D_{ij}^\circ$ dominates and the convergence slows
down. Our scheme can, on its own, not compete with the mix of
techniques presented by Bremer and Gimbutas~\cite{Brem12}.

\section{Conclusion}

This note is about promoting singularity subtraction as a helpful tool
in the discretization of singular integral operators on curved
surfaces. Singular and nearly singular kernels are expanded in series
whose terms are integrated on parametrically rectangular regions using
high-order product integration, thereby reducing the need for spatial
adaptivity and precomputed weights. A simple scheme is presented and
an application to the interior Dirichlet Laplace problem on some tori
gives around ten digit accurate results using only two expansion terms
and a modest programming- and computational effort.  Further
development, including modifications as to allow for parametrically
triangular regions, is needed before the technique may find its way
into competitive solvers.

\begin{small}

\end{small}

\end{document}